\documentclass[a4paper,11pt]{article}
\usepackage{indentfirst,latexsym,bm,amsmath,amssymb,amsthm}
\usepackage[dvips]{graphicx}

\makeatletter\@addtoreset{equation}{section} \makeatother
\newtheorem{theorem}{Theorem}[section]
\newtheorem{lemma}[theorem]{Lemma}

\newtheorem{remark}{Remark}[section]

 

\title{Blow up of Solutions to Semilinear Wave Equations with variable coefficients and boundary}
\author{Yi Zhou
\thanks{School of Mathematical Sciences, Fudan University,  Shanghai 200433, P. R. China;
 Key Laboratory of Mathematics for Nonlinear Sciences (Fudan
University), Ministry of Education of China, Shanghai 200433, P. R.
China;   Shanghai Key Laboratory for Contemporary Applied
Mathematics, School of Mathematical Sciences, Fudan University,
Shanghai 200433, P. R. China; ({\tt Email:yizhou@fudan.ac.cn})} \and
Wei Han
\thanks{School of Mathematical Sciences, Fudan University,
Shanghai 200433, P. R. China; ({\tt Email:
sh\_hanweiwei1@126.com})}.}
\date{}
\begin{document}
\maketitle
\begin{abstract}
 This paper is devoted to studying the following two initial-boundary value problems for
 semilinear wave equations with variable coefficients on exterior domain with subcritical
exponent in $n$ space dimensions:
\begin{equation} \label{0.1}
u_{tt}-\partial_{i}(a_{ij}(x)\partial_{j}u)=|u|^{p},  \   \   \   \
(x,t)\in \Omega^{c}\times (0,+\infty),  \   \     n\geq 3
\end{equation}

and
\begin{equation}\label{0.2}
u_{tt}-\partial_{i}(a_{ij}(x)\partial_{j}u)=|u_{t}|^{p},  \   \   \
\ (x,t)\in \Omega^{c}\times (0,+\infty),  \   \     n\geq 1.
\end{equation}
 where $u=u(x,t)$ is a real-valued scalar unknown function in $\Omega^{c}\times [0,+\infty)
 $, here $\Omega $ is a smooth compact obstacle in $R^{n},$  $\Omega^{c} $ is its complement,  \   $ n\geq
 3$ for \eqref{0.1}  and  $ n\geq 1$  for  \eqref{0.2},
    here  $\{a_{ij}(x) \}_{i,j=1}^{n}$ denotes  a matrix valued smooth function
    of the variable $x\in \Omega^{c}$,  which takes values in the real, symmetric,
     $n\times n$ matrices,  such that for some $C>0$,
   $$  C^{-1}|\xi|^{2}\leq a_{ij}(x)\xi_{i}\xi_{j}\leq C|\xi|^{2}, \
  \forall  \xi\in R^{n}, \    \    x\in \Omega^{c},
 $$
 here and in the sequence, a repeated sum on an index  is never indicated,  and
     $$ a_{ij}(x)=\delta_{ij},  \    \      \mbox {when} \    |x|\geq R, $$
where $\delta_{ij}$ stands for the Kronecker delta function.

  The exponents $p$ satisfies $ 1<p<p_{1}(n)$ in  \eqref{0.1},  and
  $p \leq p_{2}(n)$ in \eqref{0.2}, where $p_{1}(n)$ is  the larger root of the quadratic
 equation $ (n-1)p^{2}-(n+1)p-2=0, $
  and $p_{2}(n)=\frac{2 }{n-1}+1$, respectively.  It is well-known
  that the number $p_{1}(n) $ is the critical exponent of the
   semilinear wave equation $\eqref{0.1}$, while  $p_{2}(n) $
   is the critical exponent of  $\eqref{0.2}$.

  We will establish two blowup results  for the above two initial-boundary
   value problems, it is proved that there can be no global solutions no
matter how small the initial data are,  and also we give the
lifespan estimate of solutions for above problems.


\end{abstract}
\section{Introduction}
In this paper, we will consider the blow up of solutions of the
initial-boudary value problems for the following two semilinear wave
equations on exterior domain:

\begin{equation} \label{1.1}
\left \{
\begin{array}{lllll}
u_{tt}-\partial_{i}(a_{ij}(x)\partial_{j}u)=|u|^{p},  \   \   \   \
(x,t)\in \Omega^{c}\times (0,+\infty),  \  \      n\geq 3, \       \\
 u(0,x)=\varepsilon f(x),  \   \    \     \      u_{t}(0,x)=\varepsilon g(x),   \  \  \   x\in \Omega^{c}, \\
u(t,x)|_{\partial \Omega}=0, \   \  \mbox{  for }  \  t\geq 0,
\end{array} \right.
\end{equation}

and

\begin{equation} \label{1.2}
\left \{
\begin{array}{lllll}
u_{tt}-\partial_{i}(a_{ij}(x)\partial_{j}u)=|u_{t}|^{p},  \   \   \
\ (x,t)\in \Omega^{c}\times (0,+\infty),  \   \        n\geq 1,  \       \\
u(0,x)= \varepsilon f(x),  \   \    \     \      u_{t}(0,x)=\varepsilon g(x),   \  \  \   x\in \Omega^{c}, \\
u(t,x)|_{\partial \Omega}=0,  \   \   \mbox{  for }   \  t\geq 0,
\end{array} \right.
\end{equation}
  where  $A(x)=\{a_{ij}(x) \}_{i,j=1}^{n}$ denotes
 a matrix valued smooth function of the variable $x\in \Omega^{c}$,
 which takes values in the real, symmetric, $n\times n$ matrices,
 such that for some $C>0$,
   $$  C^{-1}|\xi|^{2}\leq a_{ij}(x)\xi_{i}\xi_{j}\leq C|\xi|^{2}, \
  \forall  \xi\in R^{n}, \    \    x\in \Omega^{c},
 $$
  here and in the sequence, a repeated sum on an index is never indicated,
 and
     $$ a_{ij}(x)=\delta_{ij},  \mbox {when} \    |x|\geq R, $$
where $\delta_{ij}$ stands for the Kronecker delta function.
  $\Omega $ is a smooth compact obstacle in $R^{n}$, $\Omega^{c} $ is its complement,
  \  $ n\geq  3$ for \eqref{1.1} and   $ n\geq 1$  for \eqref{1.2}.
Without loss of generality,  we assume that $0\in
\Omega\subset\subset B_{R}, $  where $B_{R}$ is a ball of radius $R
$ centered at the origin and  $supp\{ f, g \}\subset B_{R}$.
 We consider dimensions $n\geq 3$ and  exponents $p\in (1,
 p_{1}(n))$ for problem \eqref{1.1},  and  dimensions $n\geq 1$ and
  exponents  $ p \leq p_{2}(n)$ for problem  \eqref{1.2},
   where $p_{1}(n) $ is the larger root of the quadratic
 equation  $ (n-1)p^{2}-(n+1)p-2=0, $  and $p_{2}(n)=\frac{2
 }{n-1}+1$, respectively.   The number $p_{1}(n) $ is known as the critical exponent of the
   semilinear wave equation $\eqref{1.1}$ (see, e.g., \cite{Strauss 2}) and the number $p_{2}(n) $
   is known as the critical exponent of the
   semilinear wave equation $\eqref{1.2}$ (see, e.g., \cite{Y. Zhou 4}).
   And  we consider
 compactly supported nonnegative data $(f,g)\in
H^{1}(\Omega^{c})\times L^{2}(\Omega^{c})$  for problem \eqref{1.1}
and  $f, g\in C^{\infty}_{0}(\Omega^{c})$ for problem \eqref{1.2}.

If $a_{ij}=\delta_{ij}$,  we say problems  \eqref{1.1}, \eqref{1.2}
 are of constant coefficients. In the case of cauchy problems of
subcritical  semilinear wave equation with constant coefficients,
there is an  extensive literature which we shall review briefly,
for details, see
 \cite{ Georgiev,  Glassey 1,  Glassey 2,  John 1, JK,  Kato 1, Ta-tsien,  Schaeffer,
  Sideris, Strauss 1, Strauss 2, Takamura, Takamura 1,  Q. S. Zhang 1, Q. S. Zhang 3,   Y. Zhou 1,
 Y. Zhou 2, Y. Zhou 3, Y. Zhou
 4}.

 For the problem \eqref{1.1} with constant coefficients,  the case $n=3$ was first done by F. John
\cite{John 1} in 1979,  he showed that  when $n=3$ global solutions
always exist if $ p>p_{1}(3)=1+\sqrt{2} $ and initial data are
suitably small, and moreover, the global solutions do not exist if
$1<p<p_{1}(3)=1+\sqrt{2} $ for any nontrivial choice of $f$ and $g$.
 The number $p_{1}(3)=1+\sqrt{2} $ appears to have first arisen in
Strauss' work on low energy scattering for the nonlinear
Klein-Gordon equation \cite{Strauss 1}. This led him to conjecture
that when $n\geq 2$ global solutions of \eqref{1.1} should always
exist if initial data are sufficiently small and $p$ is greater than
a critical power $p_{1}(n)$. The conjecture was verified when $n=2$
by R. T. Glassey \cite{Glassey 2}.   In higher space dimensions, the
case $n=4$  was proved by Y. Zhou \cite{Y. Zhou 3} and  V. Georgiev,
H. Lindblad and C. Sogge \cite{Georgiev} showed that when $n\geq 4$
and $ p_{1}(n)<p\leq\frac{n+3}{n-1}$,  \eqref{1.1} has global
solutions for small initial values (see also \cite{LS} and
\cite{Tataru}). Later,  a simple proof was given by Tataru
\cite{Tataru} in the case $p>p_{1}(n) $ and $n\geq 4$. R. T. Glassey
\cite{Glassey 1} and  T. C. Sideris \cite{Sideris} showed the
blow-up result of $1<p<p_{1}(n)$ for $n=2 $  and  all $n\geq 4 $,
respectively. Sideris' proof of the blow up result is quite
delicate, using sophisticated computation involving spherical
harmonics and other special functions. His proof was simplified by
Rammaha \cite{Rammaha} and Jiao and Zhou \cite{JZ1}. In 2005, the
proof  was further simplified by Yordanov and Zhang \cite{Q. S.
Zhang 1} by using a simple test function, also, more importantly
they use their method to establish blowup phenomenon for wave
equations \eqref{1.1} with constant coefficients and a potential. On
the other hand,  for the critical case $p=p_{1}(n)$, it was shown by
Schaeffer \cite{Schaeffer} that the critical power also belongs to
the blowup case  for small data when $n=2,  3$  (see also
\cite{Takamura, Y. Zhou 1,Y. Zhou 2}).  \    B. Yordanov, and Q. S.
Zhang \cite{Q. S. Zhang
  2} and Y. Zhou\cite{Y. Zhou 5} independently  have extended   Sideris' blowup result
  to $p=p_{1}(n)$ for all $n\geq 4$ by different methods respectively.

For the problem \eqref{1.2} with constant coefficients,  the blowup
part was first proved by F. John \cite{John 3}  and  the global
existence part  was first obtained by T.C. Sideris \cite{Sideris 1}
in the case $n=3$, and both by J. Schaeffer \cite{Schaeffer 1} in
the case $n=5$.  The blow-up part in the case $n=2$ was proved by
Schaeffer \cite{Schaeffer 2} for
 $p=p_{2}(2)$. Later, R. Agemi \cite{Agemi 1} proved it for $1<p\leq
 p_{2}(2)$ by different method from \cite{Schaeffer 2}. The case
 $n=1$ is essentially due to K. Masuda \cite{Masuda} who proved the
 blowup result in the case $n=1,2,3$ and $p=2$. In higher space dimensions,
 M. A. Rammaha \cite{Rammaha} proved the blow-up part of
 $n\geq 4 $ in the case where $p=p_{2}(n)$ for odd $n$ and
 $1<p<p_{2}(n)$ for even $n$. A simple proof of blowup part was
 later given by Y. Zhou \cite{Y. Zhou 4}.

 Recently, K. Hidano et. al \cite{Hidano} has established global
 existence for problem \eqref{1.1} with $p>p_{1}(n)$ and $n=3,4$.
 For related result,  one can see Sogge and Wang's  work
 \cite{Sogge}.  However, to the best of our knowledge, there are no blowup  results
concerning initial-boundary value problems for semilinear wave
equations with variable coefficients on exterior domain.  In this
paper,  we shall establish blowup results for the initial-boundary
   value problem for subcritical values of $p$. We shall also estimate the lifespan
   $T(\varepsilon)$ for small initial data of size $\varepsilon $.
Our result is complement to the global existence result of K. Hidano
et. al \cite{Hidano}.   For the problem \eqref{1.1}, we obtain our
result by constructing
   two test functions  $\phi_{0}$ and $\psi_{1}$ (see Section 2),
   which is motivated by the work of Yordanov and Zhang
\cite{Q. S. Zhang 1}. For the problem \eqref{1.2}, we still use the
test function $\psi_{1}$ and by introducing an auxiliary function
$G_{0}(t)$ (see Section 4), we reduced the problem to a Ricatti
equation. This proof is new even in the constant coefficients case.

We are interested in showing the "blow up"  of solutions to
 problems \eqref{1.1} and \eqref{1.2}.  For that, we require
   \begin{equation}
 1<p<p_{1}(n) \  \    \mbox{for} \  \eqref{1.1},  \     \mbox{and} \  \
    p \leq p_{2}(n) \   \    \mbox{for} \    \eqref{1.2},
\end{equation}
 where $p_{1}(n) $ is the larger root of the quadratic
 equation  $ (n-1)p^{2}-(n+1)p-2=0, $  and $p_{2}(n)=\frac{2 }{n-1}+1$.
  We are also interested in estimating the time when "blow
  up" occurs. For initial data of the form
 \begin{equation}
 u(0,x)=\varepsilon f(x), \    \     \    \
 u_{t}(0,x)=\varepsilon g(x),
\end{equation}
with constant $0<\varepsilon \leq 1$, smallness can be measured
conveniently by the size of $\varepsilon $ for fixed $f$, $g$.  We
define "life span" $T(\varepsilon)$ of the solutions of \eqref{1.1}
or \eqref{1.2} to be the largest value such that solutions exist for
$x\in \Omega^{c}$,  $0\leq t< T(\varepsilon)$.

 For problem $ (1.1)$,  we consider compactly supported nonnegative data $(f,g)\in
H^{1}(\Omega^{c})\times L^{2}(\Omega^{c})$,  $n\geq 3$ and satisfy
\begin{equation}\label{1.3}
f(x)\geq 0,  \   \    g(x)\geq 0, \   a.e.,     \   \ f(x)=g(x)=0, \
 \   \mbox{ for } \   |x|>R,  \  \      \mbox { and}  \  \    f(x)\not\equiv 0.
\end{equation}
  We establish the following theorem for \eqref{1.1}:

\noindent\begin{theorem}\label{thm:1.1}   Let $(f,g)\in
H^{1}(\Omega^{c})\times L^{2}(\Omega^{c})$ and satisfy \eqref{1.3},
 $\partial \Omega  $ is smooth, and $ \Omega$ satisfies the
exterior ball conditions, space dimensions $n\geq 3$. Suppose
  that problem (1.1) has a solution $(u, u_{t})\in
C([0,T), \ H^{1}(\Omega^{c})\times L^{2}(\Omega^{c}))$ such that
 $$\emph{supp}(u, u_{t})\subset \{ (x,t): \    |x|\leq t+R    \}\cap (\Omega^{c}\times R^{+}).$$
  If $1<p<p_{1}(n)$,    then $T<\infty$, and there exists a
 positive constant $A_{1} $ which is independent of $ \varepsilon$
 such that
\begin{equation}
\begin{array}{ll}
  T(\varepsilon)&  \leq
  A_{1}\varepsilon^{-\frac{2p(p-1)}{2+(n+1)p-(n-1)p^{2} }  }.
\end{array}
\end{equation}
\end{theorem}

\begin{remark}\label{remark:1.1}
Exterior ball  condition may not be necessary, but in certain point
of our proof, we use strong maximum principle for the elliptic
equation, so this condition is needed technically.
\end{remark}

For problem \eqref{1.2},  we consider compactly supported
nonnegative data  $f, g\in C^{\infty}_{0}(\Omega^{c})$, $n\geq 1$
and satisfy
\begin{equation}\label{1.4}
 f(x)\geq 0 , \    g(x)\geq 0, \      \      \   \ f(x)=g(x)=0, \
 \   \mbox{ for}  \   |x|>R \   \   \mbox{and} \  \    g(x)  \not\equiv
 0.
\end{equation}

 Similarly, we establish the following theorem for \eqref{1.2}:
 \vskip .3cm
\noindent\begin{theorem}\label{thm:1.2}   Let  $f, g$  are smooth
functions with compact support $f, g\in C^{\infty}_{0}(\Omega^{c})$
and satisfy \eqref{1.4},   space dimensions $n\geq 1$. Suppose
  that problem (1.2) has a solution $(u, u_{t})\in
C([0,T), \ H^{1}(\Omega^{c})\times L^{q}(\Omega^{c}))$ , where $q=
\max (2, p)$ \   such that
 $$\emph{supp}(u, u_{t})\subset \{ (x,t): \    |x|\leq t+R    \}\cap (\Omega^{c}\times R^{+}).$$
 If $p\leq p_{2}(n)$,     then $T<\infty$,  moreover,  we have the following estimates for the life span
$T(\varepsilon)$ of solutions of \eqref{1.2}:
\\ (i) \    If $(n-1)(p-1)<2 $, then there exists a positive constant
 $A_{2}$ which is independent of $\varepsilon$
such that
\begin{equation}
    T(\varepsilon)\leq A_{2}
    \varepsilon^{-\frac{p-1}{1-(n-1)(p-1)/2}}.
\end{equation}
\\ (ii)\     If $(n-1)(p-1)=2$, then there exist a positive constant
$B_{2}$ which is independent of $\varepsilon$ such that
  \begin{equation}
    T(\varepsilon)\leq
    \exp({B_{2} \varepsilon^{-(p-1)}}).
\end{equation}
\end{theorem}

The rest of the paper is arranged as follows. We state several
preliminary propositions in Section 2, Section 3 is devoted to the
blowup proof  for our Theorem \ref{thm:1.1} and we prove the
Theorem \ref{thm:1.2} in Section 4.

\section{Preliminaries} \vskip .5cm

To prove the main results in this paper, we will employ the
following important ODE result:

\noindent \begin{lemma}\label{thm:2.1}  (see  \cite{Sideris}) Let
$p>1$, $a\geq 1$, and
 $(p-1)a>q-2$.  If $F\in C^{2}([0, T))$ satisfies
\\ (1) \   $F(t)\geq \delta (t+R)^{a}$,
 \\  (2) \   $\frac{d^{2}F(t)}{dt^{2}}\geq k (t+R)^{-q}[F(t)]^{p}$,
\\     with some positive constants $\delta $, $k$, and $R$,
then $F(t)$ will blow up in finite time,  $T<\infty$.  Furthermore,
we have the the following estimate for the life span $T(\delta)$ of
$F(t)$ :
   \begin{equation} \label{LIFESPAN}
   T(\delta)\leq c \delta^{ -\frac{(p-1)}{ (p-1)a-q+2}
   },
    \end{equation}
     where  $c$ is a positive constant depending on $k$ and $R$ but independent of $\delta$.
\end{lemma}

\noindent\begin{proof}\
 For the proof of blow up result part see Sideris \cite{Sideris}.
  We only prove the estimate of the life span of $F(t)$ as following:

Let us make a translation $ \tau= t\delta^{\frac{(p-1)}{ (p-1)a-q+2}
} $ and define
  $$H(\tau)=\delta^{\frac{(q-2)}{ (p-1)a-q+2}} F(t)=
  \delta^{\frac{(q-2)}{ (p-1)a-q+2}} F(\tau \delta^{\frac{-(p-1)}{ (p-1)a-q+2}
}), $$ then we have
\begin{equation}
\left\{
     \begin{array}{ll}
      H(\tau )\geq ( \delta^{\frac{(p-1)}{ (p-1)a-q+2} }  +\tau )^{a
      },      \cr\noalign{\vskip2mm}
       H''(\tau )\geq c ( \delta^{\frac{(p-1)}{ (p-1)a-q+2} }  +\tau
       )^{-q }H^{p}(\tau),
\end{array}
   \right.
\end{equation}
where $c$ is a positive constant.

So when $ \delta \leq R^{\frac{(p-1)a-q+2 }{p-1}}$,  easy
computation shows that
\begin{equation}
\left\{
     \begin{array}{ll}
      H(\tau )\geq \tau^{a},\   \      \cr\noalign{\vskip2mm}
      H''(\tau )\geq   c(R+\tau)^{-q}H^{p}(\tau).
\end{array}
   \right.
\end{equation}
     So   $ H(\tau
)$ will blow up in finite time and the life span of $F(t)$ satisfies
\eqref{LIFESPAN}.  This completes the proof.
\end{proof}

\noindent \begin{lemma} \label{thm:2.2}   There exists function
$\phi_{0}(x) \in C^{2}(\Omega^{c})$,  space dimensions $n\geq 3$,
satisfying the following boundary value problem:
\begin{equation} \label{2.1}
\left\{
     \begin{array}{ll}
     \partial_{i}(a_{ij}\partial_{j}\phi_{0}(x))=0,\ \    in \     \Omega^{c}, \   \
         \   n\geq 3,         \cr\noalign{\vskip2mm}
     \phi_{0}|_{\partial \Omega}=0,   \cr\noalign{\vskip2mm}
     |x|\rightarrow \infty,   \   \   \      \phi_{0}(x)\rightarrow
     1.
\end{array}
   \right.\ \    \    \    \      \
\end{equation}
 Moreover, $\phi_{0}(x) $ satisfies: for  $ \forall \  x\in \Omega^{c}, $
   $ 0<  \phi_{0}(x)< 1 $.
\end{lemma}

\noindent\begin{proof}\
 To solve $\phi_{0}(x)$,     let $ \tilde{\phi_{0}}$  be solution for
the following boundary value problem on exterior domain:

\begin{equation} \label{2.2}
\left\{
     \begin{array}{ll}
      \partial_{i}(a_{ij}\partial_{j}\tilde{\phi_{0}}(x))= 0,\ \    in \    \Omega^{c},   \   \
       \    \   n\geq 3,         \cr\noalign{\vskip2mm}
     \tilde{\phi_{0}}|_{\partial \Omega}=-1,   \cr\noalign{\vskip2mm}
     |x|\rightarrow \infty,   \   \   \      \tilde{\phi_{0}}(x)\rightarrow
     0,
\end{array}
   \right.
  \end{equation}
 this problem is well-posed,  it has unique solution $\tilde{\phi_{0}}(x)$,   and  by maximum principle,
   we can easily obtain  $ -1<\tilde{\phi_{0}}(x) < 0,$  for $\forall x\in \Omega^{c}, $
 then we can easily check that $\phi_{0}(x)=1+ \tilde{\phi_{0}}(x)$
satisfy the boundary value problem  \eqref{2.1},  This proves the
existence of $\phi_{0}$ in \eqref{2.1} and satisfies  $0<
\phi_{0}(x) < 1 $ for $\forall \   x\in \Omega^{c},  \   n\geq 3.$
 \    The proof is complete.
\end{proof}

Similarly, we have the following:

\noindent \begin{lemma} \label{thm:2.3}   There exists a  function $
\phi_{1}(x) \in C^{2}(\Omega^{c})$, space dimensions $n\geq 1$,
satisfying the following boundary value problem:
\begin{equation} \label{2.6}
   \left\{
     \begin{array}{ll}
     \partial_{i}(a_{ij}\partial_{j}\phi_{1}(x))= \phi_{1} ,\ \    in  \     \Omega^{c}, \   \
         \   n\geq 1,       \cr\noalign{\vskip2mm}
     \phi_{1}|_{\partial \Omega}=0,   \cr\noalign{\vskip2mm}
     |x|\rightarrow \infty,   \   \   \      \phi_{1}(x)\rightarrow
     \displaystyle\int_{S^{n-1}} e^{x\cdot \omega} d\omega.
\end{array}
   \right.
\end{equation}
 Moreover, $\phi_{1}(x)$ satisfies: there exists positive constant
 $C_{1}$, for  $ \forall x\in \Omega^{c}, $
    $0< \phi_{1}(x) \leq  C_{1}(1+|x|)^{-(n-1)/2}\cdot e^{|x|}
 $.
\end{lemma}

\noindent\begin{proof}\
 To solve $\phi_{1}(x)$,   let $ \tilde{\phi_{1}}$  be solution for
the following boundary value problem on exterior domain:
\begin{equation} \label{2.7}
\left\{
     \begin{array}{ll}
     \partial_{i}(a_{ij}\partial_{j}\tilde{\phi_{1}}(x))= \tilde{\phi_{1}}(x)-w(x), \  \    in \    \Omega^{c},
 \   \         \    n\geq 1,
        \cr\noalign{\vskip2mm}
     \tilde{\phi_{1}}|_{\partial \Omega}=-h(x)|_{\partial \Omega},     \cr\noalign{\vskip2mm}
     |x|\rightarrow \infty,   \   \   \      \tilde{\phi_{1}}(x)\rightarrow
     0,
\end{array}
   \right.
  \end{equation}
where $  h(x)= \displaystyle\int_{S^{n-1}} e^{x\cdot \omega}
d\omega, \    \   \       w(x)=
\partial_{i}((a_{ij}-\delta_{ij})\partial_{j}h)$,   since the
function $h $ satisfies $\Delta h=h$,  so by the condition of
$a_{ij}(x)$, we get $w(x)\in C_{c}^{\infty}(\Omega^{c})$,   so by
the theory of second order elliptic partial differential equation,
the problem  \eqref{2.7} is well-posed,  it has unique solution
$\tilde{\phi_{1}}(x)$,  then we can easily check that
$\phi_{1}(x)=h(x)+\tilde{\phi_{1}}(x)$  satisfies the boundary value
problem  \eqref{2.6}, this proves the existence of $\phi_{1}$ in
\eqref{2.6}.  To derive the estimate of $ \phi_{1}(x)$ in
$\Omega^{c} $, we rewrite the  boundary value problem \eqref{2.6} as
the following form:

\begin{equation}\label{2.8}
\left\{
     \begin{array}{ll}
     -\partial_{i}(a_{ij}\partial_{j}\phi_{1}(x))+ \phi_{1}(x)=0 ,\ \    in  \     \Omega^{c},
       \   \         \    n\geq 1,  \cr\noalign{\vskip2mm}
     \phi_{1}|_{\partial \Omega}=0,   \cr\noalign{\vskip2mm}
     |x|\rightarrow \infty,   \   \   \      \phi_{1}(x)\rightarrow
     h(x).
\end{array}
   \right.
\end{equation}

  So by maximum principle,   we can easily get
\begin{equation}\label{MP1}
\begin{array}{ll}
   \phi_{1}(x)> 0$, \   for $\forall  \  x\in \Omega^{c}.
   \end{array}
 \end{equation}
  Next we analyze   $\tilde{\phi_{1}}(x)$ in order to get the  estimation of $\phi_{1}(x)$,  we will prove that $\tilde{\phi_{1}}(x)$
  is bounded by some positive constant $C$, that is,  $|\tilde{\phi_{1}}(x)|\leq
  C$ for $\forall  \  x\in \Omega^{c} $.  Here and hereafter,
  we shall denote by $C$(or $c$) a positive constant in the estimates, and the meaning of
  $C$ (or $c$)
  may change from line to line.

 For this purpose, we rewrite problem  \eqref{2.7} as  follows:
\begin{equation} \label{2.10}
\left\{
     \begin{array}{ll}
    -\partial_{i}(a_{ij}\partial_{j}\tilde{\phi_{1}}(x))+ \tilde{\phi_{1}}(x)=w(x), \  \
      in \    \Omega^{c},    \   \         \    n\geq 1,    \cr\noalign{\vskip2mm}
     \tilde{\phi_{1}}|_{\partial \Omega}=-h(x)|_{\partial \Omega},     \cr\noalign{\vskip2mm}
     |x|\rightarrow \infty,   \   \   \      \tilde{\phi_{1}}(x)\rightarrow
     0.
\end{array}
   \right.
  \end{equation}

For the purpose of employing the maximum principle,  we denote
$C=\max\limits_{x\in \partial\Omega}|h(x)|+ \max\limits_{x\in
\Omega^{c}}|w(x)|>0$,  because the function $w(x)$ is compactly
supported function in $\Omega^{c}$,  so the above expression $C$ is
well defined. By the maximum principle, we can get the upper bound
of $\tilde{\phi_{1}}(x)$ as follows:

We rewrite the equation of $\tilde{\phi_{1}}(x)$ as following:
\begin{equation}\label{2.11}
\left\{
     \begin{array}{ll}
    -\partial_{i}(a_{ij}\partial_{j}(\tilde{\phi_{1}}(x)-C))+ (\tilde{\phi_{1}}(x)-C)=w(x)-C\leq 0, \  \
      in \    \Omega^{c},   \   \         \    n\geq 1,      \cr\noalign{\vskip2mm}
     (\tilde{\phi_{1}}-C)|_{\partial \Omega}=(-h(x)-C)|_{\partial \Omega}\leq 0,     \cr\noalign{\vskip2mm}
     |x|\rightarrow \infty,   \   \   \      (\tilde{\phi_{1}}(x)-C)\rightarrow -C\leq 0.
\end{array}
   \right.
  \end{equation}
 So we apply maximum principle to $(\tilde{\phi_{1}}(x)-C)$ , we can obtain for $\forall  x\in \Omega^{c}$,
  \    $\tilde{\phi_{1}}(x)-C\leq 0$,  that is,   $\tilde{\phi_{1}}(x)\leq C,$   in  $\Omega^{c}$.

  In a similar way, we can get    $-\tilde{\phi_{1}}(x)\leq C,$   in  $\Omega^{c}$.

Thus we conclude that $|\tilde{\phi_{1}}(x)|\leq
  C$  for any $x \in\Omega^{c}$. \\
  Hence we have  for  $ \forall  x \in\Omega^{c}$,
\begin{equation}\label{MP2}
\begin{array}{ll}
   \phi_{1}(x)= \tilde{\phi_{1}}(x)+h(x)\leq C+ h(x)\leq C' h(x)\leq C_{1} (1+ |x|)^{-(n-1)/2}\cdot
  e^{|x|}.
   \end{array}
 \end{equation}
This together with \eqref{MP1} implies that  $\phi_{1}(x)$ satisfies
\begin{equation}\label{MP2}
\begin{array}{ll}
     0 < \phi_{1}(x)\leq C_{1} (1+
   |x|)^{-(n-1)/2}\cdot e^{|x|},  \     \    \mbox{ in }    \Omega^{c}, \    n\geq 1.
   \end{array}
 \end{equation}
This proves Lemma \ref{thm:2.3}.
\end{proof}

In order to describe the following lemmas, we define the following
test  function
\begin{equation} \label{TF}
\psi_{1}(x, t)=\phi_{1}(x)e^{-t}, \  \     \forall      \  x\in
\Omega^{c}, \ t\geq 0.
\end{equation}
We have
 \noindent \begin{lemma} \label{thm:2.4}
  Let  $p>1$. Assume
that $\phi_{1}$ satisfy the conditions in  Lemma \ref{thm:2.3}, \
 $\psi_{1}(x, t) $ is as in \eqref{TF}.  Then for
$\forall \ t\geq 0$,
$$\begin{array}{ll}
 \displaystyle\int_{\Omega^{c}\cap   \{ |x|\leq t+R  \}} \left[
 \psi_{1}(x,t) \right]^{p/(p-1)}dx\leq C(t+R)^{n-1-(n-1)p'/2},
\end{array}$$
where $p'=p/(p-1)$ and $C$ is a positive constant.
\end{lemma}

\noindent\begin{proof}\
  Let $I(t)$ be the integral in Lemma \ref{thm:2.4},
by the property of $\phi_{1}(x)$, we have
\begin{equation}  \label{2.16}
\begin{array}{ll}
     I(t) &= \displaystyle\int_{\Omega^{c}\cap   \{ |x|\leq t+R  \}} \left[
 \psi_{1}(x,t) \right]^{p/(p-1)}dx = \displaystyle\int_{\Omega^{c}\cap   \{ |x|\leq t+R  \}} \left[
    \phi_{1}(x)e^{-t} \right]^{p/(p-1)}dx    \cr\noalign{\vskip 4mm} &
\leq   \displaystyle\int_{\Omega^{c}\cap   \{ |x|\leq t+R  \}}
\left[
 C_{1}(1+|x|)^{-(n-1)/2}\cdot e^{|x|} \right]^{p/(p-1)}\cdot e^{-tp'}dx
  \cr\noalign{\vskip 4mm} & \leq
 \displaystyle\int_{ \{ |x|\leq t+R  \}}
\left[
 C_{1}(1+|x|)^{-(n-1)/2}\cdot e^{|x|} \right]^{p/(p-1)}\cdot e^{-tp'}dx
\cr\noalign{\vskip 4mm} & =
  area(S^{n-1}) C_{1}^{p/(p-1)} \displaystyle\int_{0}^{t+R}
  (1+r)^{-(n-1)p'/2}\cdot e^{ p'r} r^{n-1} e^{-tp'} dr,
 \end{array}
\end{equation}
where $p'=p/(p-1)$  and  $S^{n-1}$ is the unit sphere in $R^{n}$. It
is sufficient to show that

\begin{equation}  \label{2.17}
\begin{array}{ll}
     I(t) & \leq
  C e^{-tp'} \displaystyle\int_{0}^{t+R}
  (1+r)^{n-1-(n-1)p'/2}\cdot e^{ p'r}  dr \leq
  C(t+R)^{n-1-(n-1)p'/2}.
 \end{array}
\end{equation}

This estimate is evident after splitting the last integral into two
parts, that is,

\begin{equation}  \label{2.18}
\begin{array}{ll}
      \displaystyle\int_{0}^{t+R}
  (1+r)^{n-1-(n-1)p'/2}\cdot e^{ p'r}  dr =
  \left[ \displaystyle\int_{0}^{(t+R)/2} + \displaystyle\int_{(t+R)/2}^{t+R} \right]
  (1+r)^{n-1-(n-1)p'/2}\cdot e^{ p'r}  dr.
 \end{array}
\end{equation}

$$ \begin{array}{ll}
      \displaystyle\int_{0}^{(t+R)/2}
  (1+r)^{n-1-(n-1)p'/2}\cdot e^{ p'r}  dr & \leq
  (1+t+R)^{q_{1}} \displaystyle\int_{0}^{(t+R)/2}  e^{ p'r} dr  \cr\noalign{\vskip 4mm} &
  =(1+t+R)^{q_{1}}\cdot  \displaystyle \frac{1}{p'}\left(  e^{ p'(t+R)/2}-1  \right)
\cr\noalign{\vskip 4mm} & \leq  (1+t+R)^{q_{1}}\cdot  \displaystyle
\frac{1}{p'} e^{ p'(t+R)/2} =  \displaystyle\frac{e^{ p'R/2}}{p'}
(1+t+R)^{q_{1}} e^{ p't/2},
  \end{array} $$
where $q_{1}=\max(0, n-1-(n-1)p'/2)$,  and

$$ \begin{array}{ll}
      \displaystyle\int_{(t+R)/2}^{t+R}
  (1+r)^{n-1-(n-1)p'/2}\cdot e^{ p'r}  dr & \leq
  2^{-q_{2}}(1+t+R)^{n-1-(n-1)p'/2} \displaystyle\int_{(t+R)/2}^{t+R}  e^{ p'r} dr  \cr\noalign{\vskip 4mm} &
  =2^{-q_{2}}(1+t+R)^{n-1-(n-1)p'/2}\cdot  \displaystyle \frac{1}{p'}\left(  e^{ p'(t+R)}-e^{ p'(t+R)/2}
    \right)
\cr\noalign{\vskip 4mm} & \leq   \displaystyle \frac{2^{-q_{2}}
e^{p'R}}{p'} \cdot (1+t+R)^{n-1-(n-1)p'/2} e^{ p't},
  \end{array} $$
where $ q_{2}=\min (0, n-1-(n-1)p'/2)$.

This proves Lemma  \ref{thm:2.4}.
\end{proof}

 \noindent \begin{lemma} \label{thm:2.5}
  Let  $p>1$. Assume
that $\phi_{0}$  and $\phi_{1}$ satisfy the conditions in Lemma
\ref{thm:2.2} and Lemma \ref{thm:2.3}, respectively,  $\psi_{1}(x,
t) $ is as in \eqref{TF},  $\partial\Omega$ and $\Omega$ satisfies
the conditions in Theorem \ref{thm:1.1}. Then for $\forall \ t\geq
0$,
\begin{equation}  \label{2.19}
\begin{array}{ll}
  \displaystyle\int_{\Omega^{c} \cap \{|x|\leq t+R \} }
[\phi_{0}(x)]^{-1/(p-1)}\cdot [\psi_{1}(x,t)]^{p/(p-1)}dx
  \leq  C(t+R)^{n-1-(n-1)p'/2},
\end{array}
\end{equation}
where $p'=p/(p-1)$ and $C$ is a positive constant.
\end{lemma}

\noindent\begin{proof}\
 To estimate the integral in Lemma \ref{thm:2.5},  we split it
 into two parts as follows
\begin{equation} \label{2.20}
\begin{array}{ll}
 & \displaystyle\int_{\Omega^{c} \cap \{|x|\leq t+R \} }
[\phi_{0}(x)]^{-1/(p-1)}\cdot
[\psi_{1}(x,t)]^{p/(p-1)}dx\cr\noalign{\vskip 4mm} & =
\displaystyle\int_{\Omega^{c} \cap  B_{R} }
[\phi_{0}(x)]^{-1/(p-1)}\cdot [\psi_{1}(x,t)]^{p/(p-1)}dx+
\displaystyle\int_{B_{R}^{c} \cap  \{|x|\leq t+R \} }
[\phi_{0}(x)]^{-1/(p-1)}\cdot [\psi_{1}(x,t)]^{p/(p-1)}dx
\cr\noalign{\vskip 4mm} & = I_{1}(t)+ I_{2}(t).
\end{array}
\end{equation}
 We will estimate  $I_{1}(t)$ and $I_{2}(t)$  separately.

First let us estimate  $I_{2}(t)$.   Since for $ \forall \    x\in
\Omega^{c}$, \   $
 0<\phi_{0}(x)< 1, $   we remark that there exists a constant  $ c\in (0,1)$,  such that when $ x\in
 B_{R}^{c} \cap  \{|x|\leq t+R \} $,
  $ \phi_{0}(x)\geq c $.   By Lemma  \ref{thm:2.4}, we have

\begin{equation}\label{2.21}
\begin{array}{ll}
  I_{2}(t) &= \displaystyle\int_{B_{R}^{c} \cap  \{|x|\leq t+R \} }
[\phi_{0}(x)]^{-1/(p-1)}\cdot [\psi_{1}(x,t)]^{p/(p-1)}dx
\cr\noalign{\vskip 4mm} &
 \leq  \displaystyle\int_{B_{R}^{c} \cap  \{|x|\leq t+R \} }
c^{-1/(p-1)}\cdot [\psi_{1}(x,t)]^{p/(p-1)}dx \cr\noalign{\vskip
4mm} &
  \leq  c^{-1/(p-1)}\displaystyle\int_{\Omega^{c} \cap   \{ |x|\leq t+R  \}}
\left[
 \psi_{1}(x,t) \right]^{p/(p-1)}dx
 \cr\noalign{\vskip
4mm} &\leq
 c^{-1/(p-1)} C(t+R)^{n-1-(n-1)p'/2},  \cr\noalign{\vskip
4mm} &= C_{2}(t+R)^{n-1-(n-1)p'/2}.
\end{array}
\end{equation}

Next we estimate $I_{1}(t)$.
 On the one hand, because of
smoothness of $\phi_{1}(x)$,   the first derivative of $\phi_{1}(x)$
is bounded in  $\Omega^{c} \cap  B_{R} $, this lead to
$\phi_{1}(x)=\phi_{1}(x)-\phi_{1}(y)\leq C_{3}|x-y|,  $  \   for
$\forall \   y\in \partial\Omega $.
   Therefore by taking the infimum on $\partial\Omega$ we have,
$$|\phi_{1}(x)|\leq C_{3} dist(x,\   \partial\Omega).$$
On the other hand, $\phi_{0}(x) $ obeys the maximum (minimum)
principle, and assumes its minimum value (zero) on $\partial\Omega
$, since $ \Omega$ satisfies exterior ball condition,  so by  \cite
[Hopf's Lemma, p. 330]{Evans}, it follows that, for any
  $y\in \partial\Omega, $  there exists an open ball $B\subset  \Omega^{c}
  $ with $y\in \partial B$,  then we have, for any
  $y\in \partial\Omega, $
\begin{equation}\label{2.22}
\begin{array}{ll}
\displaystyle\frac{\partial \phi_{0}}{\partial \nu}(y)>0,
\end{array}
\end{equation}
 where $\nu $ is the inner unit normal to $\Omega^{c}$ at $y$.
 By the compactness of $\partial\Omega$, we have, for $\forall \  y\in
 \partial\Omega$, we have
 $$ \displaystyle\frac{\partial \phi_{0}}{\partial \nu}(y)\geq C_{*}
 >0, $$
 where $C_{*}$ is a positive constant.

 For $ \forall \  x\in\Omega^{c} \cap  B_{R}$,  there exists a $y\in
 \partial\Omega$ such that  $(x-y) /\hspace{-0.1cm}/  \nu (y) $, i.e.,  $ \frac{(x-y )}{|x-y|} = \nu (y)
 $,  $\nu (y) $ is the outer unit normal to
$\partial\Omega$ at $y$.  \     So we have
  $$ \nabla \phi_{0}(y) \cdot \frac{(x-y) }{|x-y|}=\frac{\partial \phi_{0}(y) }{\partial \nu  }\geq C_{*}
 >0, $$
\begin{equation} \label{2.23}
\begin{array}{ll}
  \phi_{0}(x) &=  \phi_{0}(x)-\phi_{0}(y) =\displaystyle\int_{0}^{1}
  \nabla\phi_{0}(sx+(1-s)y)ds \cdot (x-y) \cr\noalign{\vskip 4mm} &
  = \displaystyle\int_{0}^{1} \nabla\phi_{0}(sx+(1-s)y)ds \cdot \frac{(x-y
  )}{|x-y|}\cdot |x-y|,
\end{array}
\end{equation}
by the continuity, for $\forall\  x\in\Omega^{c} \cap B_{R}$  and  $
|x-y|\ll 1 $, we know that $(sx+(1-s)y)$ is sufficiently close to
$y$, so we can guarantee that
$$ \nabla\phi_{0}(sx+(1-s)y) \cdot \frac{(x-y
  )}{|x-y|}  \geq \frac{1}{2}C_{*}>0. $$
 So there exists a positive constant  $\varepsilon_{0}>0 $ such that the above
expression holds for  $\forall\  x\in\Omega^{c} \cap B_{R}$ and $
dist(x, \partial\Omega)<\varepsilon_{0}$.

 We discuss in the following in two cases respectively: \\
One case is that  $ \mbox{  for }  x\in\Omega^{c} \cap  B_{R},
\mbox{
 and }  dist(x, \partial\Omega)<\varepsilon_{0}$,
 we have
 \begin{equation} \label{2.24}
\begin{array}{ll}
 | \phi_{0}(x)|\geq  \frac{1}{2} C_{*}|x-y|\geq  \frac{1}{2} C_{*}dist(x, \partial\Omega). \
 \end{array}
\end{equation}
The other case is that when $ x\in\Omega^{c} \cap  B_{R}$ ,  and $
dist(x, \partial\Omega)\geq \varepsilon_{0}, $  on the one hand,  by
the property of the function $\phi_{0}(x)$,
  there is a positive constant $c_{1}\in (0,1)$, such that
 $$  \phi_{0}(x)\geq c_{1}>0,$$
 on the other hand,  for $ x\in\Omega^{c} \cap  B_{R}$ ,  there
 definitely exists a positive constant $c'>0 $ such that
   $  dist(x, \partial\Omega)\leq c',  $
  so we have
\begin{equation}  \label{2.25}
\begin{array}{ll}
\displaystyle\frac{\phi_{0}(x)}{dist(x,
\partial\Omega)}\geq \displaystyle\frac{\phi_{0}(x)}{c'}
 \geq \displaystyle\frac{c_{1}}{c'}=c''>0,  \  \   \
  \mbox{for} \    x\in\Omega^{c} \cap  B_{R} , \  \mbox{and}\  \   dist(x,
\partial\Omega)\geq \varepsilon_{0},
\end{array}
\end{equation}
that is $$ \phi_{0}(x)\geq c''dist(x,
\partial\Omega),    $$
 where $c'' $ is a positive constant. \\
 So combining the above two cases,
 for $ \forall \   x\in \Omega^{c} \cap  B_{R}$, we have
   $$ \phi_{0}(x)\geq C_{**}dist(x,
\partial\Omega),       $$
 where $C_{**}$ is a positive constant. \\
 Hence, we have
\begin{equation}  \label{2.26}
\begin{array}{ll}
I_{1}(t)& =\displaystyle\int_{\Omega^{c} \cap  B_{R} }
[\phi_{0}(x)]^{-1/(p-1)}\cdot [\psi_{1}(x,t)]^{p/(p-1)}dx
\cr\noalign{\vskip 4mm} & \leq \displaystyle\int_{\Omega^{c} \cap
B_{R} } [C_{**}]^{-1/(p-1)}[dist(x,
\partial\Omega)]^{-1/(p-1)}\cdot
[\psi_{1}(x,t)]^{p/(p-1)}dx \cr\noalign{\vskip 4mm} & =
\displaystyle\int_{\Omega^{c} \cap B_{R} }
[C_{**}]^{-1/(p-1)}[dist(x,
\partial\Omega)]^{-1/(p-1)}\cdot e^{-tp'}
[\phi_{1}(x)]^{p/(p-1)}dx \cr\noalign{\vskip 4mm} & \leq
\displaystyle\int_{\Omega^{c} \cap B_{R} }
[C_{**}]^{-1/(p-1)}[dist(x,
\partial\Omega)]^{-1/(p-1)}\cdot e^{-tp'}
C_{3}^{p/(p-1)} [dist(x, \partial\Omega)]^{p/(p-1)} dx
\cr\noalign{\vskip 4mm} & =
  e^{-tp'} \displaystyle\int_{\Omega^{c} \cap B_{R} }
  [C_{**}]^{-1/(p-1)} C_{3}^{p/(p-1)} dist(x, \partial\Omega) dx
  \cr\noalign{\vskip 4mm} & =
   C e^{-tp'} \displaystyle\int_{\Omega^{c} \cap B_{R} }
   dist(x, \partial\Omega) dx\leq C_{4} e^{-tp'},
\end{array}
\end{equation}
 where $p'=p/(p-1)$.

 So we conclude that
\begin{equation}  \label{2.27}
\begin{array}{ll}
 & \displaystyle\int_{\Omega^{c} \cap \{|x|\leq t+R \} }
[\phi_{0}(x)]^{-1/(p-1)}\cdot [\psi_{1}(x,t)]^{p/(p-1)}dx
 \cr\noalign{\vskip 4mm} & =
 I_{1}(t)+ I_{2}(t) \cr\noalign{\vskip 4mm} & \leq
   C_{4} e^{-tp'}+C_{2}(t+R)^{n-1-(n-1)p'/2}
\cr\noalign{\vskip 4mm} &
 \leq  C_{5}(t+R)^{n-1-(n-1)p'/2},
\end{array}
\end{equation}
where  $C_{5}$ is a positive constant.  The proof is complete.
\end{proof}

 \noindent \begin{lemma} \label{thm:2.6}
  Let  $p>1$. Assume  that  $\phi_{1}$ satisfies the conditions in
  Lemma \ref{thm:2.3},  $\psi_{1}(x, t) $ is as in \eqref{TF}. Then for $\forall \ t\geq 0$,

\begin{equation} \label{2.28}
  \begin{array}{ll}
     \displaystyle\int_{\Omega^{c}\cap \{|x|\leq t+R \} } \psi_{1} dx
 \leq  C(t+R)^{(n-1)/2},
  \end{array}
\end{equation}
where $C$ is a positive constant.
\end{lemma}

\noindent\begin{proof}\
 We note that for $\forall \   t\geq 0 $,   $
 \psi_{1}(x,t)=e^{-t}\phi_{1}(x)$,  and since  $\mbox { for} \    \forall \   x\in  \Omega^{c}
 ,$  $ 0< \phi_{1}(x) \leq  C_{1} (1+|x|)^{-(n-1)/2} e^{|x|} $,  we
 can get that there exists a positive constant $C_{6}$ such
 that    $ 0< \phi_{1}(x) \leq  C_{6} |x|^{-(n-1)/2} e^{|x|} $  for any $x\in
 \Omega^{c}$. \\
   So  we have
\begin{equation} \label{2.29}
  \begin{array}{ll}
     \displaystyle\int_{\Omega^{c}\cap \{|x|\leq t+R \} } \psi_{1} dx
  &=  \displaystyle\int_{\Omega^{c}\cap \{|x|\leq t+R \} }  e^{-t} \phi_{1}(x)
  dx  \cr\noalign{\vskip 4mm}
     &\leq  \displaystyle\int_{\Omega^{c}\cap \{|x|\leq t+R \} }  e^{ -t}\cdot C_{6} |x|^{-(n-1)/2} e^{|x|} dx
\cr\noalign{\vskip 4mm} &\leq   \displaystyle\int_{ \{|x|\leq t+R \}
}  e^{ -t}\cdot C_{6} |x|^{-(n-1)/2} e^{|x|} dx  \cr\noalign{\vskip
4mm}
     &= C_{6} e^{ -t} \displaystyle\int_{0}^{t+R} r^{-(n-1)/2}
     e^{r}\cdot r^{n-1}dr \displaystyle \int_{S^{n-1}} d\omega= C_{7} e^{
     -t} \displaystyle\int_{0}^{t+R} e^{r}\cdot r^{(n-1)/2} dr \cr\noalign{\vskip 4mm}
     &= C_{7} e^{ -t} \left[ e^{r}r^{(n-1)/2}|_{0}^{t+R}-
     \displaystyle\int_{0}^{t+R} e^{r}(\frac{n-1}{2}) r^{(n-3)/2}dr
     \right] \cr\noalign{\vskip 4mm} &\leq
      C_{7} e^{-t}e^{t+R}(t+R)^{(n-1)/2}= C_{7} e^{R}(t+R)^{(n-1)/2}=
      C_{8}(t+R)^{(n-1)/2}.
  \end{array}
\end{equation}
This completes the proof.
\end{proof}

\section{The proof of Theorem  \ref{thm:1.1}} \vskip .5cm

Theorem \ref{thm:1.1} is a consequence
 of the lower bound and the blowup result about nonlinear
 differential inequalities in Lemma  \ref{thm:2.1}.

To outline the method, we will introduce the following functions:
\begin{equation}  \label{3.1}
\left\{
     \begin{array}{ll}
     F_{0}(t)= \displaystyle\int_{\Omega^{c}} u(x,t)\phi_{0}(x)dx,\ \     \cr\noalign{\vskip4mm}
     F_{1}(t)= \displaystyle\int_{\Omega^{c}} u(x,t)\psi_{1}(x, t)dx,\
     \  \     \                  \psi_{1}(x, t)=\phi_{1}(x)e^{-t},
\end{array}
   \right.
\end{equation}
 here $\phi_{0}(x)$ and $\phi_{1}(x)$ are as in Lemma
 \ref{thm:2.2} and Lemma \ref{thm:2.3}.   The assuptions on $u$ imply that $F_{0}(t)$ and $F_{1}(t)$ are
 well-defined $C^{2}$-functions for all $t$. By a standard
 procedure, we derive a nonlinear differential inequality for
 $F_{0}(t)$. We also derive a linear differential inequality for
 $F_{1}(t)$ and combine these to obtain a polynomial lower bound on
 $F_{0}(t)$ as $t\rightarrow \infty$.

To this end, we first establish the following lemma:
 \noindent
\begin{lemma} \label{thm:3.1} Let $(f, g)$ satisfy \eqref{1.3}.
Suppose that problem \eqref{1.1} has a solution $(u, u_{t})\in
C([0,T), \ H^{1}(\Omega^{c})\times L^{2}(\Omega^{c})),
 $ such that
 $$supp(u, u_{t})\subset \{ (x,t): \    |x|\leq t+R    \}\cap (\Omega^{c}\times R^{+}).$$  Then
 for all $t\geq 0$,
  $$F_{1}(t)\geq \frac{1}{2}(1-e^{-2t})\varepsilon \displaystyle\int_{\Omega^{c}} [f(x)+g(x)]\phi_{1}(x)dx+
   e^{-2t} \varepsilon \displaystyle\int_{\Omega^{c}} f(x)\phi_{1}(x)dx\geq \varepsilon c_{0}>0.
   $$
\end{lemma}

\noindent\begin{proof}\
  We multiply \eqref{1.1} by the test function $\psi_{1}\in
C^{2}(\Omega^{c}\times R)$  and integrate over $\Omega^{c}\times
[0,t]$,  then we use integration by parts and Lemma \ref{thm:2.3}.\\
 First,
$$\begin{array}{ll}
     \displaystyle\int_{0}^{t}\int_{\Omega^{c}}
     \psi_{1}(\partial_{i}(a_{ij}(x)\partial_{j}u)- u_{tt}+|u|^{p})dxd\tau=0. \cr\noalign{\vskip 4mm}
\end{array}$$
 By the expression $\psi_{1}(x, t)=\phi_{1}(x)e^{-t} $  and Lemma
 \ref{thm:2.3}, we have
$$\begin{array}{ll}
     \displaystyle\int_{0}^{t}\int_{\Omega^{c}} \psi_{1}\partial_{i}(a_{ij}(x)\partial_{j}u) dxd\tau
&=\displaystyle\int_{0}^{t}\left[ \int_{\partial \Omega}
\psi_{1}a_{ij}(x)\partial_{j}u\cdot n_{i}  dS-
\displaystyle\int_{\Omega^{c}}
(a_{ij}(x)\partial_{i}\psi_{1})\partial_{j}u
 dx \right]d\tau \cr\noalign{\vskip 4mm} &=
-\displaystyle\int_{0}^{t}\left[ \int_{\partial \Omega}
  a_{ij}(x)\partial_{i}\psi_{1}\cdot u\cdot
  n_{j}dS-\displaystyle\int_{\Omega^{c}} \partial_{j}(a_{ij}(x)\partial_{i}\psi_{1})u
  dx \right] d\tau \cr\noalign{\vskip 4mm}
  &=\displaystyle\int_{0}^{t} \displaystyle\int_{\Omega^{c}}
  \psi_{1}u dxd\tau,
 \cr\noalign{\vskip 4mm}
  \end{array}$$
 by the expression of $\psi_{1}(x, t) $, we get $(\psi_{1})_{t}=-\psi_{1}, \   \    (\psi_{1})_{tt}=\psi_{1}
 $. So we have
$$\begin{array}{ll}
     &\displaystyle\int_{0}^{t}\int_{\Omega^{c}} \psi_{1} u_{tt}
     dxd\tau =\displaystyle\int_{0}^{t}\int_{\Omega^{c}}\left[ \partial_{\tau} (
\psi_{1}u_{\tau} )-(\psi_{1})_{\tau}u_{\tau} \right] dxd\tau
\cr\noalign{\vskip 4mm} &= \displaystyle\int_{\Omega^{c}}
\psi_{1}u_{\tau}dx |_{\tau=t} -\displaystyle\int_{\Omega^{c}}
\psi_{1}u_{\tau}dx |_{\tau=0}+
\displaystyle\int_{0}^{t}\int_{\Omega^{c}} \psi_{1}u_{\tau}dxd\tau
\cr\noalign{\vskip 4mm} & =\displaystyle\int_{\Omega^{c}}
\psi_{1}u_{\tau}dx |_{\tau=t} -\displaystyle\int_{\Omega^{c}}
\psi_{1}u_{\tau}dx |_{\tau=0}+
\displaystyle\int_{0}^{t}\int_{\Omega^{c}} \left[ \partial_{\tau}
(\psi_{1}u)-(\psi_{1})_{\tau}u  \right] dxd\tau
 \cr\noalign{\vskip 4mm} &
 =\displaystyle\int_{\Omega^{c}} \psi_{1}u_{\tau}dx |_{\tau=t}
-\displaystyle\int_{\Omega^{c}} \psi_{1}u_{\tau}dx |_{\tau=0}+
\displaystyle\int_{\Omega^{c}}
\psi_{1}udx|_{\tau=t}-\displaystyle\int_{\Omega^{c}}
\psi_{1}udx|_{\tau=0}+ \displaystyle\int_{0}^{t}\int_{\Omega^{c}}
\psi_{1}u dxd\tau
 \cr\noalign{\vskip 4mm} &
= \displaystyle\int_{\Omega^{c}} \left( \psi_{1}u_{t}+u\psi_{1}
\right) dx-\varepsilon\displaystyle\int_{\Omega^{c}}
\phi_{1}(x)g(x)dx- \varepsilon\displaystyle\int_{\Omega^{c}}
\phi_{1}(x)f(x)dx+ \displaystyle\int_{0}^{t}\int_{\Omega^{c}}
\psi_{1} u dxd\tau.
\end{array}$$

Combining the above equalities, we have
$$\begin{array}{ll}
     &\displaystyle\int_{0}^{t}\int_{\Omega^{c}}  \psi_{1}|u|^{p}
     dxd\tau= \displaystyle\int_{\Omega^{c}}\left(  \psi_{1} u_{t} +
     \psi_{1} u \right) dx-\varepsilon
     \displaystyle\int_{\Omega^{c}}\phi_{1}(x) [f(x)+g(x)]dx.
\end{array}$$

We notice that
$$\begin{array}{ll}
     \displaystyle\int_{\Omega^{c}}  \left(  \psi_{1}u_{t} + \psi_{1}u
     \right) dx&= \displaystyle\frac{d}{dt}\displaystyle\int_{\Omega^{c}} (
     \psi_{1} u )dx- \displaystyle\int_{\Omega^{c}} (\psi_{1})_{t}u dx +
   \displaystyle\int_{\Omega^{c}} \psi_{1}u dx \cr\noalign{\vskip 4mm}&
 = \displaystyle\frac{d}{dt}\displaystyle\int_{\Omega^{c}} ( \psi_{1} u
 )dx+ 2 \displaystyle\int_{\Omega^{c}} \psi_{1}u
 dx\cr\noalign{\vskip 4mm}& =\displaystyle\frac{dF_{1}(t)}{dt}+2F_{1}(t).
\end{array}$$
So  by $\psi_{1}>0$, we have
$$\begin{array}{ll}
     \displaystyle\frac{dF_{1}(t)}{dt}+2F_{1}(t)&=
 \displaystyle\int_{0}^{t}\int_{\Omega^{c}}
     |u|^{p}\psi_{1}(x,\tau)dxd\tau + \varepsilon\displaystyle\int_{\Omega^{c}} \phi_{1}(x)\left[
f(x)+g(x) \right]dx \cr\noalign{\vskip 2mm}& \geq
\varepsilon\displaystyle\int_{\Omega^{c}} \left[ f(x)+g(x)
\right]\phi_{1}(x)dx.
\end{array}$$
Multiplying the above expression by $e^{2t}$,   we obtain
$$\begin{array}{ll}
  \displaystyle\frac{d(e^{2t}F_{1}(t))}{dt}\geq e^{2t} \varepsilon\displaystyle\int_{\Omega^{c}} \left[ f(x)+g(x)
\right]\phi_{1}(x)dx,
 \cr\noalign{\vskip
2mm}
\end{array}$$
and integrating the above differential inequality over $[0,t]$, we
get

$$\begin{array}{ll}
 e^{2t}F_{1}(t)-F_{1}(0)\geq \displaystyle\frac{1}{2}(e^{2t}-1)\varepsilon \displaystyle\int_{\Omega^{c}} \left[ f(x)+g(x)
\right]\phi_{1}(x)dx.
 \cr\noalign{\vskip
2mm}
\end{array}$$
   Observing $F_{1}(0)=\displaystyle\int_{\Omega^{c}}u(x,0)\psi_{1}(x, 0)dx
   = \varepsilon\displaystyle\int_{\Omega^{c}}f(x)\phi_{1}(x)dx $.
So, by the property of the function $f(x)$ and $\phi_{1}(x)$,  we
arrive at
$$\begin{array}{ll}
 F_{1}(t)\geq  \displaystyle\frac{1}{2}(1-e^{-2t})\varepsilon\displaystyle\int_{\Omega^{c}} \left[ f(x)+g(x)
\right]\phi_{1}(x)dx+
e^{-2t}\varepsilon\displaystyle\int_{\Omega^{c}}
f(x)\phi_{1}(x)dx\geq \varepsilon c_{0}>0 .
 \cr\noalign{\vskip
2mm}
\end{array}$$
Thus we obtain the lower bound in Lemma \ref{thm:3.1}.
\end{proof}

   Next We shall show that $F_{0}(t)$ satisfies the
 differential inequalities in Lemma  \ref{thm:2.1} for suitable $a,\
 q$.    For this purpose,  we  multiply  \eqref{1.1} by $\phi_{0}$  and integrate over $\Omega^{c}$.
  We note that for a fixed $t$, $u(\cdot , t)\in H^{1}_{0}(D_{t})$ where
 $D_{t}$ is the support of $u(\cdot , t)$. Hence we can use integration by parts and Lemma
 \ref{thm:2.2}. \\
   First,
$$\begin{array}{ll}
     \displaystyle\int_{\Omega^{c}} \left[ \phi_{0} \partial_{i}(a_{ij}(x)\partial_{j}u)- \phi_{0}
u_{tt}+|u|^{p}\phi_{0} \right] dx=0. \cr\noalign{\vskip 4mm}
\end{array}$$

Since \begin{equation}
\begin{array}{ll} \label{3.2}
       \displaystyle\int_{\Omega^{c}} \phi_{0}
       \partial_{i}(a_{ij}(x)\partial_{j}u) dx  &= \displaystyle\int_{\partial \Omega}
   \phi_{0}a_{ij}(x)\partial_{j}u\cdot n_{i} dS- \displaystyle\int_{\Omega^{c}}
    \partial_{i}\phi_{0}a_{ij}(x)\partial_{j}u dx   \cr\noalign{\vskip 4mm}
   &=-\left( \displaystyle\int_{\partial \Omega}  a_{ij}(x)
   \partial_{i}\phi_{0}u\cdot n_{j} dS- \displaystyle\int_{\Omega^{c}}
    \partial_{j}(a_{ij}(x)\partial_{i}\phi_{0})u dx  \right) \cr\noalign{\vskip 4mm}
   &= \displaystyle\int_{\Omega^{c}} \partial_{j}(a_{ij}(x)\partial_{i}\phi_{0})u dx
   =0.
 \end{array}
\end{equation}
So  we get  $$ \frac{d^{2}F_{0}(t)}{dt^{2}} =
\displaystyle\int_{\Omega^{c}} |u(x,t)|^{p}\phi_{0}(x) dx.  $$

Estimating the right side of the above equality by the
H$\ddot{o}$lder inequality, we have
$$\begin{array}{ll}
    & \left| \displaystyle\int_{\Omega^{c}} u(x,t)\phi_{0}(x) dx
    \right| \cr\noalign{\vskip 4mm}
     &= \left| \displaystyle\int_{\Omega^{c} \cap \{|x|\leq t+R \}} u(x,t)[\phi_{0}(x)]^{1/p} [\phi_{0}(x)]^{(p-1)/p} dx
    \right| \cr\noalign{\vskip 4mm}
&\leq   \left( \displaystyle\int_{\Omega^{c} \cap \{|x|\leq t+R \} }
|u(x,t)[\phi_{0}(x)]^{1/p}|^{p} dx\right)^{1/p} \cdot \left(
\displaystyle\int_{\Omega^{c} \cap \{|x|\leq t+R \} }
|[\phi_{0}(x)]^{(p-1)/p}|^{p'} dx\right)^{1/p'} \cr\noalign{\vskip
2mm}
\end{array}$$
 where $p'= p/(p-1)$,
this implies that
$$\begin{array}{ll}
    \left| \displaystyle\int_{\Omega^{c}} u(x,t)\phi_{0}(x) dx
    \right|^{p} &\leq  \left( \displaystyle
    \int_{ \{|x|\leq t+R \} \cap \Omega^{c}} |u(x,t)|^{p}\phi_{0}(x)dx\right)
     \left(  \displaystyle\int_{\{|x|\leq t+R \} \cap \Omega^{c}} \phi_{0}(x) dx
     \right)^{p-1}  \cr\noalign{\vskip 4mm}
       & \leq  \left( \displaystyle\int_{\Omega^{c}} |u(x,t)|^{p}\phi_{0}(x)dx\right)
     \left(  \displaystyle\int_{ \{|x|\leq t+R \} \cap \Omega^{c}} \phi_{0}(x) dx
     \right)^{p-1}. \cr\noalign{\vskip 2mm}
\end{array}$$
 So we have
$$\begin{array}{ll}
    \displaystyle\int_{\Omega^{c}} |u(x,t)|^{p}\phi_{0}(x)dx &\geq
    \displaystyle\frac{ \left| \displaystyle\int_{\Omega^{c}} u(x,t)\phi_{0}(x) dx
    \right|^{p} }{ \left(  \displaystyle\int_{\{|x|\leq t+R \} \cap \Omega^{c}} \phi_{0}(x) dx
     \right)^{p-1}}.   \cr\noalign{\vskip 2mm}
\end{array}$$

By Lemma 2.2,  we have
$$\begin{array}{ll}
  \displaystyle\int_{\{|x|\leq t+R \} \cap \Omega^{c}} \phi_{0}(x) dx \leq
 \displaystyle\int_{\{|x|\leq t+R \} } 1 dx
\leq   Vol \{x: |x|\leq
  t+R \} =Vol(\mathbf{B}^{n})(t+R)^{n}.
       \cr\noalign{\vskip 2mm}
\end{array}$$

Therefore
$$\begin{array}{ll}
  \displaystyle\int_{\Omega^{c}} |u(x,t)|^{p}\phi_{0}(x)dx \geq
  \displaystyle\frac{ \left| \int_{\Omega^{c}} u(x,t)\phi_{0}(x) dx
    \right|^{p} }{  [ Vol(\mathbf{B}^{n})(t+R)^{n}]^{p-1}}=
  \frac{ \left| F_{0}(t) \right|^{p} }{  [Vol(\mathbf{B}^{n})]^{p-1}\cdot
  (t+R)^{n(p-1)}}. \cr\noalign{\vskip 2mm}
\end{array}$$

Thus
\begin{equation}\label{3.3}
\begin{array}{ll}
  \displaystyle\frac{d^{2}F_{0}(t)}{dt^{2}} \geq k
  (t+R)^{-n(p-1)}\cdot \left| F_{0}(t) \right|^{p},
  \cr\noalign{\vskip 2mm}
\end{array}
\end{equation}
 where $k= [Vol(\mathbf{B}^{n})]^{-(p-1)}>0.$
So $F_{0}$ satisfies the differential inequality (2) in Lemma
\ref{thm:2.1}. To show that $F_{0}$ admits the lower bound (1) in
Lemma \ref{thm:2.1}, we relate $d^{2}F_{0}(t)/dt^{2}$ to $F_{1}$
using again \eqref{1.1} and the H$\ddot{o}$lder inequality.

Since
$$\begin{array}{ll}
   & \left| \displaystyle\int_{\Omega^{c}} u(x,t)\psi_{1}(x,t) dx
    \right| \cr\noalign{\vskip 4mm} &= \left| \displaystyle\int_{\Omega^{c} \cap \{|x|\leq t+R \} } u(x,t)[\phi_{0}(x)]^{1/p}\cdot
     [\phi_{0}(x)]^{-1/p} \cdot \psi_{1}(x,t)  dx
    \right| \cr\noalign{\vskip 4mm}
&\leq \left( \displaystyle\int_{\Omega^{c} \cap \{|x|\leq t+R \} }
|u(x,t)|^{p} \cdot \phi_{0}(x) dx\right)^{1/p} \cdot \left(
\displaystyle\int_{\Omega^{c} \cap \{|x|\leq t+R \} }
|[\phi_{0}(x)]^{-1/p}\cdot \psi_{1}(x,t)|^{p'} dx\right)^{1/p'}
\cr\noalign{\vskip 2mm}
  &\leq \left( \displaystyle\int_{\Omega^{c}}
|u(x,t)|^{p}\cdot \phi_{0}(x) dx\right)^{1/p} \cdot \left(
\displaystyle\int_{\Omega^{c} \cap \{|x|\leq t+R \} }
[\phi_{0}(x)]^{-1/(p-1)}\cdot
[\psi_{1}(x,t)]^{p/(p-1)}dx\right)^{(p-1)/p},  \cr\noalign{\vskip
2mm}
\end{array}$$
 where $p'= p/(p-1)$, this implies that
$$\begin{array}{ll}
 & \left| \displaystyle\int_{\Omega^{c}} u(x,t)\psi_{1}(x,t) dx
    \right|^{p} \cr\noalign{\vskip 4mm} &  \leq \left( \displaystyle\int_{\Omega^{c}}
|u(x,t)|^{p}\cdot \phi_{0}(x) dx\right)\cdot \left(
\displaystyle\int_{\Omega^{c} \cap \{|x|\leq t+R \} }
[\phi_{0}(x)]^{-1/(p-1)}\cdot
[\psi_{1}(x,t)]^{p/(p-1)}dx\right)^{p-1}.
  \cr\noalign{\vskip 2mm}
\end{array}$$

By  \eqref{3.1}, the above becomes
$$\begin{array}{ll}
    \displaystyle\frac{d^{2}F_{0}(t)}{dt^{2}} &= \displaystyle\int_{\Omega^{c}}
|u(x,t)|^{p}\phi_{0}(x) dx  \cr\noalign{\vskip 4mm} &
 \geq   \displaystyle\frac{\left| \displaystyle\int_{\Omega^{c}} u(x,t)\psi_{1}(x,t) dx
    \right|^{p} }{\left(
\displaystyle\int_{\Omega^{c} \cap \{|x|\leq t+R \} }
[\phi_{0}(x)]^{-1/(p-1)}\cdot
[\psi_{1}(x,t)]^{p/(p-1)}dx\right)^{p-1}}  \cr\noalign{\vskip 4mm} &
= \displaystyle\frac{\left| F_{1}(t) \right|^{p} }{\left(
\displaystyle\int_{\Omega^{c} \cap \{|x|\leq t+R \} }
[\phi_{0}(x)]^{-1/(p-1)}\cdot
[\psi_{1}(x,t)]^{p/(p-1)}dx\right)^{p-1}}.  \cr\noalign{\vskip 4mm}
&
\end{array}$$
In the following, we will estimate the numerator and denominator,
respectively, and provide a lower bound on  $d^{2}F_{0}/dt^{2}$.

By the Lemma \ref{thm:3.1}, we have
\begin{equation}  \label{3.4}
\left| F_{1}(t) \right|^{p} \geq \varepsilon^{p} (c_{0})^{p}>0.
  \end{equation}

Also, by the Lemma \ref{thm:2.5} we know that
\begin{equation}  \label{3.5}
\begin{array}{ll}
  \displaystyle\int_{\Omega^{c} \cap \{|x|\leq t+R \} }
[\phi_{0}(x)]^{-1/(p-1)}\cdot [\psi_{1}(x,t)]^{p/(p-1)}dx
  \leq  C_{5}(t+R)^{n-1-(n-1)p'/2},
\end{array}
\end{equation}
where $p'=p/(p-1)$ and $C_{5}$ is a positive constant.

 So by combining \eqref{3.4} and \eqref{3.5}, we obtain
$$   \begin{array}{ll}
    \displaystyle\frac{d^{2}F_{0}(t)}{dt^{2}} &\geq
 \displaystyle\frac{\left| F_{1}(t) \right|^{p} }{\left(
\displaystyle\int_{\Omega^{c} \cap \{|x|\leq t+R \} }
[\phi_{0}(x)]^{-1/(p-1)}\cdot
[\psi_{1}(x,t)]^{p/(p-1)}dx\right)^{p-1}}.  \cr\noalign{\vskip 4mm}
& \geq  \displaystyle\frac{ \varepsilon^{p} c_{0}^{p} }{\left[
 C_{5}(t+R)^{n-1-(n-1)p'/2} \right]^{p-1}}  \cr\noalign{\vskip 4mm}
& \geq  L (t+R)^{-(n-1)(p/2-1)},
\end{array}$$
where  $L=\varepsilon^{p}c_{0}^{p}C_{5}^{-(p-1)}>0. $   Integrating
twice, we have the final estimate

 $$ F_{0}(t)\geq  \delta (t+R)^{n+1-(n-1)p/2}+\displaystyle\frac{dF_{0}(0)}{dt}t+F_{0}(0),  $$
 with  constant $$ \delta=\displaystyle\frac{ L }{ [n-\frac{1}{2}(n-1)p+1 ][n-\frac{1}{2}(n-1)p]  }=
  \displaystyle\frac{\varepsilon^{p}c_{0}^{p}C_{5}^{-(p-1)} }{ [n-\frac{1}{2}(n-1)p+1 ][n-\frac{1}{2}(n-1)p] }
  >0.
  $$
     When $1<p<p_{1}(n)$, it is easy to check that  $ n+1-(n-1)p/2>1.
   $ Hence the following estimate is valid when $t$ is sufficiently large:
\begin{equation} \label{3.6}
\begin{array}{ll}
  F_{0}(t)\geq  \frac{1}{2} \delta (t+R)^{n+1-(n-1)p/2}.
\end{array}
\end{equation}

 Estimates  \eqref{3.3} together with \eqref{3.6} and Lemma \ref{thm:2.1} with parameters
   $$ a\equiv  n+1-(n-1)p/2, \   \mbox{and}  \    q\equiv n(p-1) $$
   imply Theorem \ref{thm:1.1} for all exponents $p$ such that
    $$(p-1)(n+1-(n-1)p/2)> n(p-1)-2 \   \      \mbox{and} \  \      p>1.    $$
    It is easy to see that the solution set is $p\in (1,p_{1}(n))$,
    so by Lemma  \ref{thm:2.1}, all solutions of problem \eqref{1.1}
   with nontrivial nonnegative initial values must blow up in finite time.

   Also, recall from Lemma \ref{thm:2.1}, we have the following estimate for the life
    span $T(\varepsilon)$ of solutions of \eqref{1.1} as follows:

\begin{equation}
\begin{array}{ll}
  T(\varepsilon)&  \leq  c (\frac{1}{2} \delta )^{ -\frac{(p-1)}{ (p-1)a-q+2}
  } \cr\noalign{\vskip2mm} &
  = A_{1} \left( \varepsilon^{p} \right)^{ -\frac{(p-1)}{ (p-1)(
  n+1-(n-1)p/2)-n(p-1)+2} } \cr\noalign{\vskip 4mm} &
  =A_{1}\varepsilon^{-\frac{p(p-1)}{ (p-1)(
  1-(n-1)p/2)+2}  }\cr\noalign{\vskip 4mm} &
   =A_{1}\varepsilon^{-\frac{2p(p-1)}{ 2+(n+1)p-(n-1)p^{2}} },
\end{array}
\end{equation}
 where $A_{1} $ is a positive constant which is independent of $ \varepsilon$.
    The proof of Theorem \ref{thm:1.1} is complete.

\section{The proof of Theorem  \ref{thm:1.2}} \vskip .4cm

   By the expression   $ \psi_{1}(x,t)=e^{-t}\phi_{1}(x)\geq 0$,     \
  we have  $ (\psi_{1})_{t} =-\psi_{1} $,   and
   $\partial_{i}(a_{ij}\partial_{j}\psi_{1}(x, t))= \psi_{1}(x, t), \  \    in  \
   \Omega^{c}\times (0,+\infty)$.  So $\psi_{1}$ satisfies
\begin{equation} \label{4.1}
   \left\{
     \begin{array}{ll}
     \partial_{i}(a_{ij}\partial_{j}\psi_{1}(x,t))= \psi_{1}(x,t) ,\ \    in  \  \     \Omega^{c}\times (0,+\infty),     \cr\noalign{\vskip2mm}
     \psi_{1}|_{\partial \Omega\times (0,+\infty)}=0,   \cr\noalign{\vskip2mm}
     |x|\rightarrow \infty,   \   \   \      \psi_{1}(x,t)\rightarrow
     e^{-t}\displaystyle\int_{S^{n-1}} e^{x\cdot \omega} d\omega,  \   \   \mbox{for }   \  t\geq 0.
\end{array}
   \right.
\end{equation}

We multiply \eqref{1.2} by  function $\psi_{1}$,   and integrate
over $\Omega^{c}$,  then we use  integration by parts and Lemma
\ref{thm:2.3}.

First,
$$\begin{array}{ll}
     \displaystyle\int_{\Omega^{c}} \psi_{1}(u_{tt}-\partial_{i}(a_{ij}(x)\partial_{j}u))dx =
     \displaystyle\int_{\Omega^{c}} \psi_{1}|u_{t}|^{p}dx. \cr\noalign{\vskip 4mm}
\end{array}$$
 Note that for a fixed $t$, $u(\cdot , t)\in H^{1}_{0}(D_{t})$, where
 $D_{t}$ is the support of $u(\cdot , t)$. Hence by integration by parts and Lemma  \ref{thm:2.3},  we
 have
$$\begin{array}{ll}
     \displaystyle\int_{\Omega^{c}} \psi_{1}\partial_{i}(a_{ij}(x)\partial_{j}u ) dx
&= \displaystyle\int_{\Omega^{c}}(
\partial_{i}[\psi_{1}a_{ij}(x)\partial_{j}u ]-\partial_{i}\psi_{1}a_{ij}(x)\partial_{j}u
) dx  \cr\noalign{\vskip 4mm} &= \displaystyle\int_{\partial \Omega}
\psi_{1}a_{ij}(x)\partial_{j}u\cdot n_{i}  dS-
\displaystyle\int_{\Omega^{c}}
(a_{ij}(x)\partial_{i}\psi_{1})\cdot\partial_{j}u
 dx \cr\noalign{\vskip 4mm} &=
-\displaystyle \int_{\partial \Omega}
a_{ij}(x)\partial_{i}\psi_{1}\cdot u\cdot
  n_{j}dS + \displaystyle\int_{\Omega^{c}} \partial_{j}(a_{ij}(x)\partial_{i}\psi_{1})u
  dx  \cr\noalign{\vskip 4mm}
  &=\displaystyle\int_{\Omega^{c}}
  \partial_{j}(a_{ij}(x)\partial_{i}\psi_{1})\cdot u dx=\displaystyle\int_{\Omega^{c}}
  \psi_{1}\cdot u dx.
 \cr\noalign{\vskip 4mm}
  \end{array}$$
Combining the above two identities, we conclude
\begin{equation} \label{4.2}
\begin{array}{ll}
     \displaystyle\int_{\Omega^{c}} \psi_{1}u_{tt}-
     \displaystyle\int_{\Omega^{c}} \psi_{1}\cdot u dx=
     \displaystyle\int_{\Omega^{c}} \psi_{1} \cdot |u_{t}|^{p}dx. \cr\noalign{\vskip 4mm}
\end{array}
\end{equation}

Notice that
\begin{equation}\label{4.3}
     \displaystyle\frac{d }{dt} \displaystyle\int_{\Omega^{c}}
     \psi_{1}u_{t}dx= \displaystyle\int_{\Omega^{c}} ( \psi_{1}\cdot u_{tt}-u_{t}\psi_{1})
     dx,
\end{equation}

\begin{equation} \label{4.4}
     \displaystyle\frac{d }{dt} \displaystyle\int_{\Omega^{c}}
     (\psi_{1}u) dx= \displaystyle\int_{\Omega^{c}}  \left[ (\psi_{1})_{t}\cdot u+u_{t}\cdot
     \psi_{1} \right]dx= \displaystyle\int_{\Omega^{c}}
      \left[ \psi_{1}\cdot u_{t}-u \psi_{1}\right]dx.
\end{equation}
Adding up the above two expression, we obtain the following
\begin{equation} \label{4.5}
     \displaystyle\frac{d }{dt} \displaystyle\int_{\Omega^{c}}
     \left(\psi_{1}u_{t}+ \psi_{1}u \right) dx= \displaystyle\int_{\Omega^{c}}
  \left( \psi_{1}\cdot u_{tt}-u\cdot
     \psi_{1} \right) dx=  \displaystyle\int_{\Omega^{c}} \psi_{1} \cdot |u_{t}|^{p}dx.
\end{equation}
So we have
\begin{equation} \label{4.6}
\begin{array}{ll}
     \displaystyle\int_{\Omega^{c}}
     \left(\psi_{1}u_{t}+ \psi_{1}u \right) dx &= \displaystyle\int_{\Omega^{c}}
      \left(\psi_{1}u_{t}+ \psi_{1}u \right) dx|_{t=0}+
\displaystyle\int_{0}^{t}\int_{\Omega^{c}}
   \psi_{1}\cdot |u_{\tau}|^{p}dxd \tau \cr\noalign{\vskip 4mm} &
= \displaystyle\int_{\Omega^{c}}
      \varepsilon \phi_{1}(x)\left[f(x)+ g(x) \right] dx+ \displaystyle\int_{0}^{t}\int_{\Omega^{c}}
   \psi_{1}\cdot |u_{\tau}|^{p}dxd \tau \cr\noalign{\vskip 4mm} &
\geq   \displaystyle\int_{\Omega^{c}}
     \varepsilon \phi_{1}(x)g(x) dx+ \displaystyle\int_{0}^{t}\int_{\Omega^{c}}
   \psi_{1}\cdot |u_{\tau}|^{p}dxd \tau.
 \end{array}
\end{equation}
Adding two expressions  \eqref{4.2} and \eqref{4.6}, we have
\begin{equation} \label{4.7}
\begin{array}{ll}
     \displaystyle\int_{\Omega^{c}}
     \left(\psi_{1}u_{tt}+ \psi_{1}u_{t} \right) dx &\geq \displaystyle\int_{\Omega^{c}}
      \psi_{1}\cdot |u_{t}|^{p}dx +
\displaystyle\int_{0}^{t}\int_{\Omega^{c}}
   \psi_{1}\cdot |u_{\tau}|^{p}dxd \tau +\varepsilon \displaystyle\int_{\Omega^{c}}
      \phi_{1}(x)g(x) dx.
 \end{array}
\end{equation}

 Also, we know that
\begin{equation} \label{4.8}
\begin{array}{ll}
     \displaystyle\frac{d }{dt} \displaystyle\int_{\Omega^{c}}
     \psi_{1}u_{t} dx + 2 \displaystyle\int_{\Omega^{c}}
     \psi_{1}\cdot u_{t}dx&= \displaystyle\int_{\Omega^{c}}
  \left[ \psi_{1}u_{tt}+u_{t}(\psi_{1})_{t}+2\psi_{1}u_{t} \right]
  dx \cr\noalign{\vskip 4mm} &= \displaystyle\int_{\Omega^{c}} \left( \psi_{1}u_{tt}+ \psi_{1}u_{t}
  \right)dx.
 \end{array}
\end{equation}
 So we have
\begin{equation} \label{4.9}
\begin{array}{ll}
     \displaystyle\frac{d }{dt} \displaystyle\int_{\Omega^{c}}
     \psi_{1}u_{t} dx + 2 \displaystyle\int_{\Omega^{c}}
     \psi_{1}\cdot u_{t}dx&\geq  \displaystyle\int_{\Omega^{c}}
      \psi_{1}\cdot |u_{t}|^{p}dx +
\displaystyle\int_{0}^{t}\int_{\Omega^{c}}
   \psi_{1}\cdot |u_{\tau}|^{p}dxd \tau +\varepsilon \displaystyle\int_{\Omega^{c}}
      \phi_{1}(x)g(x) dx.
 \end{array}
\end{equation}

 To show the blowup property, we define the following  auxiliary function
\begin{equation} \label{4.10}
     G_{0}(t)=\displaystyle\int_{\Omega^{c}} \psi_{1}u_{t} dx-\displaystyle
     \frac{1}{2}\displaystyle\int_{0}^{t}\int_{\Omega^{c}}
   \psi_{1}\cdot |u_{\tau}|^{p}dxd \tau-\displaystyle
     \frac{\varepsilon}{2}\displaystyle\int_{\Omega^{c}}
      \phi_{1}(x)g(x) dx.
\end{equation}
  We  note that, when $t=0$,
   $$G_{0}(0)=\varepsilon\displaystyle\int_{\Omega^{c}}
      \phi_{1}(x)g(x) dx-\displaystyle
     \frac{\varepsilon}{2}\displaystyle\int_{\Omega^{c}}
      \phi_{1}(x)g(x) dx=\displaystyle
     \frac{\varepsilon}{2}\displaystyle\int_{\Omega^{c}}
      \phi_{1}(x)g(x) dx\geq 0, $$
  and we have
  \begin{equation}\label{4.11}
      \displaystyle\frac{d }{dt} G_{0}(t)=\displaystyle\int_{\Omega^{c}}
      \left( \psi_{1}u_{tt}-u_{t}\psi_{1} \right) dx-\displaystyle
     \frac{1}{2}\displaystyle\int_{\Omega^{c}}
   \psi_{1}\cdot |u_{t}|^{p}dx.
\end{equation}
  Hence,  we conclude that
 \begin{equation} \label{4.12}
  \begin{array}{ll}
     & \displaystyle\frac{d }{dt} G_{0}(t)+2G_{0}(t)  \cr\noalign{\vskip 4mm} &=\displaystyle\int_{\Omega^{c}}
      \left( \psi_{1}u_{tt}+u_{t}\psi_{1} \right) dx-\displaystyle
     \frac{1}{2}\displaystyle\int_{\Omega^{c}}
   \psi_{1}\cdot |u_{t}|^{p}dx- \displaystyle\int_{0}^{t}\int_{\Omega^{c}}
   \psi_{1}\cdot |u_{\tau}|^{p}dxd \tau - \varepsilon\displaystyle\int_{\Omega^{c}}
      \phi_{1}(x)g(x) dx  \cr\noalign{\vskip 4mm} &
 \geq   \displaystyle\int_{\Omega^{c}} \psi_{1}\cdot |u_{t}|^{p}dx +
 \displaystyle\int_{0}^{t}\int_{\Omega^{c}}\psi_{1}\cdot |u_{\tau}|^{p}dxd
 \tau + \varepsilon \displaystyle\int_{\Omega^{c}} \phi_{1}(x)g(x) dx
        - \displaystyle
     \frac{1}{2}\displaystyle\int_{\Omega^{c}} \psi_{1}
     |u_{t}|^{p}dx
       \cr\noalign{\vskip 3mm} &
         \    \     \     \      - \displaystyle\int_{0}^{t}
     \int_{\Omega^{c}}\psi_{1}\cdot |u_{\tau}|^{p}dxd \tau - \varepsilon \displaystyle\int_{\Omega^{c}}
      \phi_{1}(x)g(x) dx  \cr\noalign{\vskip 4mm} &
     =    \displaystyle \frac{1}{2}
     \displaystyle\int_{\Omega^{c}} \psi_{1}\cdot |u_{t}|^{p}dx\geq
     0.
 \end{array}
\end{equation}
   Multiplying the above differential inequality by $e^{2t}$,  we get the following expression
   $$\displaystyle\frac{d }{dt} \left( e^{2t}G_{0}(t) \right)\geq 0. $$
    So  for $ \forall \  t \geq 0,$  we have  $e^{2t}G_{0}(t)\geq G_{0}(0)
    $,   that is  $G_{0}(t)\geq e^{-2t}G_{0}(0)\geq 0. $  \\
  By  \eqref{4.10}, we have  for $ \forall t \geq 0,$
\begin{equation} \label{4.13}
    \displaystyle\int_{\Omega^{c}} \psi_{1}u_{t} dx  \geq \displaystyle
     \frac{1}{2}\displaystyle\int_{0}^{t}\int_{\Omega^{c}}
   \psi_{1}\cdot |u_{t}|^{p}dxd \tau +\displaystyle
     \frac{\varepsilon}{2}\displaystyle\int_{\Omega^{c}}
      \phi_{1}(x)g(x) dx.
\end{equation}

 Let  \begin{equation} \label{4.14}
     F(t)=\displaystyle \frac{1}{2}\displaystyle\int_{0}^{t}\int_{\Omega^{c}}
   \psi_{1}\cdot |u_{t}|^{p}dxd \tau +\displaystyle
     \frac{\varepsilon}{2}\displaystyle\int_{\Omega^{c}}
      \phi_{1}(x)g(x) dx,  \    \      t\geq 0.
\end{equation}
Then we have
   \begin{equation} \label{4.15}
      \displaystyle\int_{\Omega^{c}} \psi_{1}u_{t} dx \geq  F(t),
   \   \mbox{for} \  \    \forall  \  t \geq 0.
   \end{equation}

   Next we only need to prove that $F(t)$ blow up.\\
 From the expression of $F(t)$,  we get  for $ \forall \   t \geq 0,$ $F(t)\geq 0
 $,  and  $F'(t)= \displaystyle \frac{1}{2}\displaystyle\int_{\Omega^{c}}
   \psi_{1}\cdot |u_{t}|^{p}dx$.
  Estimating the right side of $F'(t)$ by the Holder inequality, we have

$$\begin{array}{ll}
    \displaystyle\int_{\Omega^{c}} |u_{t}(x,t)|^{p}\psi_{1}dx &\geq
    \displaystyle\frac{ \left| \displaystyle\int_{\Omega^{c}} u_{t}(x,t)\psi_{1} dx
    \right|^{p} }{  \left(  \displaystyle\int_{\Omega^{c}\cap \{|x|\leq t+R \}} \psi_{1} dx
     \right)^{p-1}},    \cr\noalign{\vskip 4mm}
\end{array}$$

 By Lemma \ref{thm:2.6},  we know that for $\forall \   t\geq 0 $,

\begin{equation} \label{4.16}
  \begin{array}{ll}
     \displaystyle\int_{\Omega^{c}\cap \{|x|\leq t+R \} } \psi_{1} dx
 \leq  C_{8}(t+R)^{(n-1)/2},
  \end{array}
\end{equation}
where $C_{8}$ is a positive constant. \\
 Therefore we conclude that
 \begin{equation} \label{4.17}
 \begin{array}{ll}
 F'(t)&= \displaystyle \frac{1}{2}\displaystyle\int_{\Omega^{c}}
   \psi_{1}\cdot |u_{t}|^{p}dx \geq
   \displaystyle \frac{1}{2} \displaystyle\frac{ \left| \displaystyle\int_{\Omega^{c}}
u_{t}(x,t)\psi_{1} dx \right|^{p} }{  \left(
\displaystyle\int_{\Omega^{c}\cap \{|x|\leq t+R \}} \psi_{1} dx
     \right)^{p-1}}  \cr\noalign{\vskip 4mm} & \geq
  C_{9} \displaystyle\frac{ \left| \displaystyle\int_{\Omega^{c}}
u_{t}(x,t)\psi_{1} dx \right|^{p} }{  (t+R )^{(n-1)(p-1)/2}}\geq
C_{9} \displaystyle\frac{ \left| F(t) \right|^{p} }{  (t+R
)^{(n-1)(p-1)/2}} .
 \end{array}
\end{equation}

 By the property of Ricatti equation, we know that  when $(n-1)(p-1)/2\leq
 1$,  the solution of the initial-boundary value problem \eqref{1.2} blow
 up.

 In detail, let $$ M= \displaystyle \frac{1}{2}\displaystyle\int_{\Omega^{c}}
 \phi_{1}(x) g(x)dx.  $$
 Then  $F(t)$ satisfies the following problem

\begin{equation} \label{4.18}
   \left\{
     \begin{array}{ll}
    F'(t) \geq  C_{9} \displaystyle\frac{ \left| F(t)  \right|^{p} }{  (t+R )^{(n-1)(p-1)/2}},
    \cr\noalign{\vskip2mm}
    F(0)=M \varepsilon.
\end{array}
   \right.
\end{equation}
 Now we introduce a function $v(t)$  satisfying the following
 Ricatti equation
\begin{equation} \label{4.19}
   \left\{
     \begin{array}{ll}
    v'(t) =  C_{9} \displaystyle\frac{ \left| v(t)  \right|^{p} }{  (t+R )^{(n-1)(p-1)/2}},
    \cr\noalign{\vskip2mm}
    v(0)=M \varepsilon.
\end{array}
   \right.
\end{equation}
  So the life span of $F$  is less than that of $v$ which will be
  the upper bound of $T(\varepsilon )$.

 Thus, in the case $ (n-1)(p-1)<2, $  integrating \eqref{4.19}, we get
\begin{equation}  \label{4.20}
     \begin{array}{ll}
    v(t) = \left[ (M\varepsilon )^{-(p-1)}  +C'R^{1-(n-1)(p-1)/2} -C'(t+R)^{1-(n-1)(p-1)/2}
          \right]^{-\frac{1}{p-1} },
\end{array}
\end{equation}
 where $$C'=\frac{C_{9}(p-1) }{1-(n-1)(p-1)/2}.   $$
  Thus
     $$  T(\varepsilon)\leq A_{2} \varepsilon^{-\frac{p-1}{1-(n-1)(p-1)/2}},$$
  where  $A_{2}$ is a positive constant  which is independent of
  $\varepsilon$.

     When $(n-1)(p-1)=2$, integrating \eqref{4.19}, we get
\begin{equation} \label{4.21}
       \begin{array}{ll}
    v(t) = \left[ (M\varepsilon )^{-(p-1)}-C'' \ln \left( \frac{t+R}{R }\right)
        \right]^{-\frac{1}{p-1} },
\end{array}
  \end{equation}
 where $$C''=C_{9}(p-1). $$
   $$   T(\varepsilon)\leq
    \exp({B_{2} \varepsilon^{-(p-1)}}),   $$
where  $B_{2}$ is a positive constant which is independent of
  $\varepsilon$.  This ends the proof of Theorem \ref{thm:1.2}.

\section*{Acknowledgments.}
This work is supported by the National Natural Science Foundation of
China (10728101), the 973 Project of the Ministry of the Science and
Technology of China, the Doctoral Foundation of the Ministry of
Education of China and the `111' Project (B08018) and SGST
09DZ2272900.

\end{document}